\documentclass[a4paper]{amsart}
\usepackage[T1]{fontenc}
\usepackage[utf8]{inputenc}
\usepackage[active]{srcltx}
\usepackage{amsthm}
\usepackage{amssymb}
\usepackage{esint}

\makeatletter

\newcommand{\noun}[1]{\textsc{#1}}

\theoremstyle{plain}
\newtheorem{thm}{\protect\theoremname}
  \theoremstyle{remark}
  \newtheorem*{rem*}{\protect\remarkname}

\makeatother

  \providecommand{\remarkname}{Remark}
\providecommand{\theoremname}{Theorem}

\begin{document}

\title{A Note on Boole Polynomials}

\author{Dae San Kim and Taekyun Kim}
\begin{abstract}
Boole polynomials play an important role in
the area of number theory, algebra and umbral calculus. In this paper,
we investigate some properties of Boole polynomials and consider
Witt-type formulas for the Boole numbers and polynomials. Finally,
we derive some new identities of those polynomials from the Witt-type
formulas which are related to Euler polynomials.
\end{abstract}
\maketitle

\section{Introduction}

\global\long\def\Zp{\mathbb{Z}_{p}}
\global\long\def\Qp{\mathbb{Q}_{p}}
\global\long\def\Cp{\mathbb{C}_{p}}

Let $p$ be a fixed prime number. Throughout this paper, $\Zp,$ $\Qp$
and $\Cp$ will denote the ring of $p$-adic integers, the field of
$p$-adic numbers and the completion of algebraic closure of $\Qp$.
The $p$-adic norm $\left|\cdot\right|_{p}$ is normalized as $\left|p\right|_{p}=\frac{1}{p}$.
Let $C\left(\Zp\right)$ be the space of continuous functions on $\Zp$.
For $f\in C\left(\Zp\right)$, the fermionic $p$-adic integral on
$\Zp$ is defined by Kim to be 
\begin{equation}
I\left(f\right)=\int_{\Zp}f\left(x\right)d\mu\left(x\right)=\lim_{N\rightarrow\infty}\sum_{x=0}^{p^{N}-1}f\left(x\right)\left(-1\right)^{x},\label{eq:1}
\end{equation}
(see \cite{key-9}).

From (\ref{eq:1}), we can derive 
\begin{equation}
I\left(f_{1}\right)=-I\left(f\right)+2f\left(0\right),\label{eq:2}
\end{equation}
(see {[}1,2,4-8,12,14-16{]}). Here $f_{1}(x)=f(x+1)$.

For $k\in\mathbb{N}$, the Euler polynomials of order $k$
are defined by the generating function to be 
\begin{equation}
\left(\frac{2}{e^{t}+1}\right)^{k}e^{xt}=\sum_{n=0}^{\infty}E_{n}^{\left(k\right)}\left(x\right)\frac{t^{n}}{n!},\label{eq:3}
\end{equation}
(see \cite{key-4,key-6,key-9,key-16}). 

When $x=0$, $E_{n}^{\left(k\right)}=E_{n}^{\left(k\right)}\left(0\right)$
are called the Euler numbers of order $k$.

In particular, for $k=1$, $E_{n}\left(x\right)=E_{n}^{\left(1\right)}\left(x\right)$
are called the Euler polynomials. 

The Stirling numbers of the first kind are defined by 
\begin{equation}
\left(x\right)_{n}=x\left(x-1\right)\cdots\left(x-n+1\right)=\sum_{l=0}^{n}S_{1}\left(n,l\right)x^{l},\label{eq:4}
\end{equation}
(see \cite{key-13}). 

As is well known, the Stirling numbers of the second kind are given
by the generating function to be 
\begin{equation}
\left(e^{t}-1\right)^{n}=n!\sum_{l=n}^{\infty}S_{2}\left(l,n\right)\frac{t^{l}}{l!},\label{eq:5}
\end{equation}
(see \cite{key-13}). 

The Boole polynomials are defined by the generating function to
be 
\begin{equation}
\sum_{n=0}^{\infty}Bl_{n}\left(x|\lambda\right)\frac{t^{n}}{n!}=\frac{1}{\left(1+\left(1+t\right)^{\lambda}\right)}\left(1+t\right)^{x},\label{eq:6}
\end{equation}
(see \cite{key-3,key-10,key-11,key-13}).

When $\lambda=1$, $2Bl_{n}\left(x|1\right)=Ch_{n}\left(x\right)$
are Changhee polynomials which are defined by 
\[
\sum_{n=0}^{\infty}Ch_{n}\left(x\right)\frac{t^{n}}{n!}=\frac{2}{t+2}\left(1+t\right)^{x},
\]
(see \cite{key-7}).

In this paper, we investigate some properties of Boole polynomials
and consider Witt-type formulas for the Boole numbers and polynomials.
Finally, we derive some new identities of those polynomials from the
Witt-type formulas which are related to Euler polynomials.

\section{A Note on Boole Polynomials}

In this section, we assume that $t\in\Cp$ with $\left|t\right|_{p}<p^{-\frac{1}{p-1}}$
and $\lambda\in\Zp$. Let us take $f\left(x\right)=\left(1+t\right)^{\lambda x}.$
From (\ref{eq:2}), we have 
\begin{eqnarray}
\int_{\Zp}\left(1+t\right)^{x+\lambda y}d\mu_{-1}\left(y\right) & = & \frac{2}{1+\left(1+t\right)^{\lambda}}\left(1+t\right)^{x}\label{eq:7}\\
 & = & \sum_{n=0}^{\infty}2Bl_{n}\left(x|\lambda\right)\frac{t^{n}}{n!}.\nonumber 
\end{eqnarray}
It is easy to show that 
\begin{equation}
\int_{\Zp}\left(1+t\right)^{x+\lambda y}d\mu_{-1}\left(y\right)=\sum_{n=0}^{\infty}\int_{\Zp}\left(x+y\lambda\right)_{n}d\mu_{-1}\left(y\right)\frac{t^{n}}{n!}.\label{eq:8}
\end{equation}

Therefore, by (\ref{eq:7}) and (\ref{eq:8}), we obtain the following
theorem. 
\begin{thm}
\label{thm:1}For $n\ge0$, we have 
\[
\int_{\Zp}\left(x+y\lambda\right)_{n}d\mu_{-1}\left(y\right)=2Bl_{n}\left(x|\lambda\right).
\]

\end{thm}
By (\ref{eq:6}), we get 
\begin{equation}
\sum_{n=0}^{\infty}2Bl_{n}\left(x|\lambda\right)\frac{\left(e^{t}-1\right)^{n}}{n!}=\frac{2}{e^{\lambda t}+1}e^{xt}=\sum_{m=0}^{\infty}E_{m}\left(\frac{x}{\lambda}\right)\lambda^{m}\frac{t^{m}}{m!},\label{eq:9}
\end{equation}
and 
\begin{eqnarray}
\sum_{n=0}^{\infty}2Bl_{n}\left(x|\lambda\right)\frac{1}{n!}\left(e^{t}-1\right)^{n} & = & \sum_{n=0}^{\infty}2Bl_{n}\left(x|\lambda\right)\frac{1}{n!}n!\sum_{m=n}^{\infty}S_{2}\left(m,n\right)\frac{t^{m}}{m!}\label{eq:10}\\
 & = & \sum_{m=0}^{\infty}\left(2\sum_{n=0}^{\infty}Bl_{n}\left(x|\lambda\right)S_{2}\left(m,n\right)\right)\frac{t^{m}}{m!}.\nonumber 
\end{eqnarray}

Therefore, by (\ref{eq:9}) and (\ref{eq:10}), we obtain the following
theorem.
\begin{thm}
For $m\ge0$, 
\[
\sum_{n=0}^{m}Bl_{n}\left(\lambda|x\right)S_{2}\left(m,n\right)=\frac{1}{2}E_{m}\left(\frac{x}{\lambda}\right)\lambda^{m}.
\]
\end{thm}
\begin{rem*}
From Theorem \ref{thm:1}, we note that 
\begin{eqnarray*}
2Bl_{n}\left(x|\lambda\right) & = & \sum_{l=0}^{n}S_{1}\left(n,l\right)\int_{\Zp}\left(x+y\lambda\right)^{l}d\mu_{-1}\left(y\right)\\
 & = & \sum_{l=0}^{\infty}S_{1}\left(n,l\right)\lambda^{l}E_{l}\left(\frac{x}{\lambda}\right).
\end{eqnarray*}

\end{rem*}
$\:$

Let us consider the Boole polynomials of order $k$$\left(\in\mathbb{N}\right)$
as follows : 
\begin{equation}
2^{k}Bl_{n}^{\left(k\right)}\left(x|\lambda\right)=\int_{\Zp}\cdots\int_{\Zp}\left(\lambda x_{1}+\cdots+\lambda x_{k}+x\right)_{n}d\mu_{-1}\left(x_{1}\right)\cdots d\mu_{-1}\left(x_{k}\right).\label{eq:11}
\end{equation}
Thus, by (\ref{eq:11}), we get 
\begin{equation}
2^{k}Bl_{n}^{\left(k\right)}\left(x|\lambda\right)=\sum_{l=0}^{n}S_{1}\left(n,l\right)\lambda^{l}E_{l}^{\left(k\right)}\left(\frac{x}{\lambda}\right).\label{eq:12}
\end{equation}

From (\ref{eq:11}), we can derive the generating function of $Bl_{n}^{\left(k\right)}\left(x\right)$
as follows : 
\begin{align}
 & 2^{k}\sum_{n=0}^{\infty}Bl_{n}^{\left(k\right)}\left(x|\lambda\right)\frac{t^{n}}{n!}\label{eq:13}\\
= & \int_{\Zp}\cdots\int_{\Zp}\left(1+t\right)^{\lambda x_{1}+\cdots+\lambda x_{k}+x}d\mu_{-1}\left(x_{1}\right)\cdots d\mu_{-1}\left(x_{k}\right)\nonumber \\
= & \left(\frac{2}{1+\left(1+t\right)^{\lambda}}\right)^{k}\left(1+t\right)^{x}.\nonumber 
\end{align}
By replacing $t$ by $e^{t}-1$, we get 
\begin{eqnarray}
\sum_{n=0}^{\infty}2^{k}Bl_{n}^{\left(k\right)}\left(x|\lambda\right)\frac{1}{n!}\left(e^{t}-1\right)^{n} & = & \left(\frac{2}{e^{\lambda t}+1}\right)^{k}e^{xt}\label{eq:14}\\
 & = & \sum_{m=0}^{\infty}\lambda^{m}E_{m}^{\left(k\right)}\left(\frac{x}{\lambda}\right)\frac{t^{m}}{m!},\nonumber 
\end{eqnarray}
and 
\begin{equation}
\sum_{n=0}^{\infty}2^{k}Bl_{n}^{\left(k\right)}\left(x|\lambda\right)\frac{\left(e^{t}-1\right)^{n}}{n!}=\sum_{m=0}^{\infty}\left(\sum_{n=0}^{m}2^{k}Bl_{n}^{\left(k\right)}\left(x|\lambda\right)S_{2}\left(m,n\right)\right)\frac{t^{m}}{m!}.\label{eq:15}
\end{equation}

Therefore, by (\ref{eq:14}) and (\ref{eq:15}), we obtain the following
theorem.
\begin{thm}
For $m\ge0$, we have 
\[
\sum_{n=0}^{m}Bl_{n}^{\left(k\right)}\left(x|\lambda\right)S_{2}\left(m,n\right)=\frac{\lambda^{m}}{2^{k}}E_{m}^{\left(k\right)}\left(\frac{x}{\lambda}\right).
\]

\end{thm}
For $n\ge0$, the rising factorial sequence is defined by 
\begin{eqnarray}
x^{\left(n\right)} & = & x\left(x+1\right)\cdots\left(x+n-1\right)=\left(-1\right)^{n}\left(-x\right)_{n}\label{eq:16}\\
 & = & \sum_{l=0}^{n}\begin{bmatrix}n\\
l
\end{bmatrix}x^{l},\nonumber 
\end{eqnarray}
where $\begin{bmatrix}n\\
l
\end{bmatrix}=\left(-1\right)^{n-l}S_{1}\left(n,l\right).$

Now, we define the Boole polynomials of the second kind as follows:
\global\long\def\bls{\hat{B}l}
\begin{equation}
\bls_{n}\left(x|\lambda\right)=\frac{1}{2}\int_{\Zp}\left(-\lambda y+x\right)_{n}d\mu_{-1}\left(y\right),\quad\left(n\ge0\right).\label{eq:17}
\end{equation}

Then, by (\ref{eq:17}), we get 
\begin{eqnarray}
\bls_{n}\left(x|\lambda\right) & = & \frac{1}{2}\sum_{l=0}^{n}S_{1}\left(n,l\right)\left(-1\right)^{l}\lambda^{l}\int_{\Zp}\left(-\frac{x}{\lambda}+y\right)^{l}d\mu_{-1}\left(y\right)\label{eq:18}\\
 & = & \frac{1}{2}\sum_{l=0}^{n}S_{1}\left(n,l\right)\left(-1\right)^{l}\lambda^{l}E_{l}\left(-\frac{x}{\lambda}\right).\nonumber 
\end{eqnarray}

When $x=0$, $\bls_{n}\left(\lambda\right)=\bls_{n}\left(0|\lambda\right)$
are called the Boole numbers of the second kind. 

The generating function of $\bls_{n}\left(x|\lambda\right)$ is given
by 
\begin{eqnarray}
\sum_{n=0}^{\infty}\bls_{n}\left(x|\lambda\right)\frac{t^{n}}{n!} & = & \frac{1}{2}\int_{\Zp}\left(1+t\right)^{-\lambda y+x}d\mu_{-1}\left(y\right)\label{eq:19}\\
 & = & \frac{\left(1+t\right)^{\lambda}}{1+\left(1+t\right)^{\lambda}}\left(1+t\right)^{x}.\nonumber 
\end{eqnarray}

By replacing $t$ by $e^{t}-1$, we get 
\begin{equation}
\sum_{n=0}^{\infty}\bls_{n}\left(x|\lambda\right)\frac{1}{n!}\left(e^{t}-1\right)^{n}=\sum_{m=0}^{\infty}\frac{\lambda^{m}}{2}E_{m}\left(\frac{\lambda+x}{\lambda}\right)\frac{t^{m}}{m!},\label{eq:20}
\end{equation}
and 
\begin{equation}
\sum_{n=0}^{\infty}\bls_{n}\left(x|\lambda\right)\frac{1}{n!}\left(e^{t}-1\right)^{n}=\sum_{m=0}^{\infty}\left(\sum_{n=0}^{m}\bls_{n}\left(x|\lambda\right)S_{2}\left(m,n\right)\right)\frac{t^{m}}{m!}.\label{eq:21}
\end{equation}

Therefore, by (\ref{eq:20}) and (\ref{eq:21}), we obtain the following
theorem.
\begin{thm}
For $m\ge0$, we have 
\[
\frac{\lambda^{m}}{2}E_{m}\left(\frac{\lambda+x}{\lambda}\right)=\sum_{n=0}^{m}\bls_{n}\left(x|\lambda\right)S_{2}\left(m,n\right),
\]
and 
\[
\bls_{m}\left(x|\lambda\right)=\sum_{l=0}^{m}S_{1}\left(m,l\right)\left(-1\right)^{l}\frac{\lambda^{l}}{2}E_{l}\left(-\frac{x}{\lambda}\right).
\]

\end{thm}
For $k\in\mathbb{N}$, let us consider the Boole polynomials of the
second kind with order $k$ as follows : 
\begin{align}
 & \bls_{n}^{\left(k\right)}\left(x|\lambda\right)\label{eq:22}\\
= & \frac{1}{2^{k}}\int_{\Zp}\cdots\int_{\Zp}\left(-\left(\lambda x_{1}+\cdots+\lambda x_{k}\right)+x\right)_{n}d\mu_{-1}\left(x_{1}\right)\cdots d\mu_{-1}\left(x_{k}\right).\nonumber 
\end{align}

Thus, by (\ref{eq:22}), we get

\begin{equation}
2^{k}\bls_{n}^{\left(k\right)}\left(x|\lambda\right)=\sum_{l=0}^{n}S_{1}\left(n,l\right)\lambda^{l}\left(-1\right)^{l}E_{l}^{\left(k\right)}\left(-\frac{x}{\lambda}\right).\label{eq:23}
\end{equation}

The generating function of $\bls_{n}^{\left(k\right)}\left(x|\lambda\right)$
is given by 
\begin{align}
 & \sum_{n=0}^{\infty}\bls_{n}^{\left(k\right)}\left(x|\lambda\right)\frac{t^{n}}{n!}\label{eq:24}\\
= & \frac{1}{2^{k}}\int_{\Zp}\cdots\int_{\Zp}\left(1+t\right)^{-\left(\lambda x_{1}+\cdots+\lambda x_{k}\right)+x}d\mu_{-1}\left(x_{1}\right)\cdots d\mu_{-1}\left(x_{k}\right)\nonumber \\
= & \left(\frac{\left(1+t\right)^{\lambda}}{1+\left(1+t\right)^{\lambda}}\right)^{k}\left(1+t\right)^{x}.\nonumber 
\end{align}

By replacing $t$ by $e^{t}-1$, we get 
\begin{eqnarray}
\sum_{n=0}^{\infty}\bls_{n}^{\left(k\right)}\left(x|\lambda\right)\frac{1}{n!}\left(e^{t}-1\right)^{n} & = & \frac{1}{2^{k}}\left(\frac{2}{e^{\lambda t}+1}\right)^{k}e^{\left(\lambda k+x\right)t}\label{eq:25}\\
 & = & \sum_{m=0}^{\infty}\frac{\lambda^{m}}{2^{k}}E_{m}^{\left(k\right)}\left(k+\frac{x}{\lambda}\right)\frac{t^{m}}{m!}\nonumber 
\end{eqnarray}
and 
\begin{equation}
\sum_{n=0}^{\infty}\bls_{n}^{\left(k\right)}\left(x|\lambda\right)\frac{1}{n!}\left(e^{t}-1\right)^{n}=\sum_{m=0}^{\infty}\left(\sum_{n=0}^{m}\bls_{n}^{\left(k\right)}\left(x|\lambda\right)S_{2}\left(m,n\right)\right)\frac{t^{m}}{m!}.\label{eq:26}
\end{equation}

Therefore, by (\ref{eq:25}) and (\ref{eq:26}), we obtain the following
theorem.
\begin{thm}
For $m\ge0$, we have 
\[
\bls_{m}^{\left(k\right)}\left(x|\lambda\right)=\sum_{l=0}^{m}S_{1}\left(m,l\right)\left(-1\right)^{l}\frac{\lambda^{l}}{2^{k}}E_{l}^{\left(k\right)}\left(-\frac{x}{\lambda}\right)
\]
and
\[ 
\frac{\lambda^{m}}{2^{k}}E_{m}^{\left(k\right)}\left(k+\frac{x}{\lambda}\right)=\sum_{n=0}^{m}\bls_{n}^{\left(k\right)}\left(x|\lambda\right)S_{2}\left(m,n\right).
\]

\end{thm}
Now, we observe that 
\begin{eqnarray}
\left(-1\right)^{n}\frac{Bl_{n}\left(x|\lambda\right)}{n!} & = & \left(-1\right)^{n}\int_{\Zp}\dbinom{x+y\lambda}{n}d\mu_{-1}\left(y\right)\label{eq:27}\\
 & = & \int_{\Zp}\dbinom{-y\lambda-x+n-1}{n}d\mu_{-1}\left(y\right)\nonumber \\
 & = & \sum_{m=0}^{n}\dbinom{n-1}{n-m}\int_{\Zp}\dbinom{-y\lambda-x}{m}d\mu_{-1}\left(y\right)\nonumber \\
 & = & \sum_{m=1}^{n}\frac{\tbinom{n-1}{m-1}}{m!}m!\int_{\Zp}\dbinom{-y\lambda-x}{m}d\mu_{-1}\left(y\right)\nonumber \\
 & = & \sum_{m=1}^{n}\dbinom{n-1}{m-1}\frac{\bls_{m}\left(-x|\lambda\right)}{m!},\nonumber 
\end{eqnarray}
and 
\begin{eqnarray}
\left(-1\right)^{n}\frac{\bls_{n}\left(x|\lambda\right)}{n!} & = & \sum_{m=0}^{n}\dbinom{n-1}{m-1}\int_{\Zp}\dbinom{-x+y\lambda}{m}d\mu_{-1}\left(y\right)\label{eq:28}\\
 & = & \sum_{m=1}^{n}\dbinom{n-1}{m-1}\frac{\bls_{m}\left(-x|\lambda\right)}{m!}.\nonumber 
\end{eqnarray}

Therefore, by (\ref{eq:27}) and (\ref{eq:28}), we obtain the following
theorem.
\begin{thm}
For $n\ge1$, we have 
\[
\left(-1\right)^{n}\frac{Bl_{n}\left(x|\lambda\right)}{n!}=\sum_{m=1}^{n}\dbinom{n-1}{m-1}\frac{\bls_{m}\left(-x|\lambda\right)}{m!}
\]
and 
\[
\left(-1\right)^{n}\frac{\bls_{n}\left(x|\lambda\right)}{n!}=\sum_{m=1}^{n}\dbinom{n-1}{m-1}\frac{Bl_{m}\left(-x|\lambda\right)}{n!}.
\]
\end{thm}

$\,$

\noindent \noun{Department of Mathematics, Sogang University, Seoul
121-742, Republic of Korea}

\noindent \emph{E-mail}\noun{ }\emph{address : }\texttt{dskim@sogang.ac.kr}

\noun{$\,$}

\noindent \noun{Department of Mathematics, Kwangwoon University, Seoul
139-701, Republic of Korea}

\noindent \emph{E-mail}\noun{ }\emph{address : }\texttt{tkkim@kw.ac.kr}
\end{document}